\documentclass[11pt,leqno]{article}
\usepackage{latexsym,amssymb}
\usepackage{amsmath}

\newtheorem{Theorem}{\bf Theorem}[section]
\newtheorem{Lemma}{\bf Lemma}[section]
\newtheorem{Proposition}{\bf Proposition}[section]
\newtheorem{Corollary}{\bf Corollary}[section]
\newtheorem{Remark}{\bf Remark}[section]
\newtheorem{Example}{\bf Example}[section]
\newtheorem{Definition}{\bf Definition}[section]

\newenvironment{theorem}{\begin{Theorem}$\!\!\!$}{\end{Theorem}}
\newenvironment{lemma}{\begin{Lemma}$\!\!\!$}{\end{Lemma}}

\newenvironment{remark}{\begin{Remark}$\!\!\!$}{\end{Remark}}

\numberwithin{equation}{section}

\begin{document}

\title{When does the heat equation have a solution\\
 with a sequence of similar level sets?}
\author{
Tatsuki Kawakami
\footnote{Department of Mathematical Sciences, Osaka Prefecture University,
Sakai 599-8531, Japan
(e-mail address : kawakami@ms.osakafu-u.ac.jp)}
\,\,\,
and
\,\,\,
Shigeru Sakaguchi
\footnote{Research Center for Pure and Applied Mathematics, Graduate School of Information Sciences, Tohoku University,
Sendai 980-8579, Japan
(e-mail address : sigersak@m.tohoku.ac.jp)}
}
\date{}
\maketitle
%    General info
%\subjclass[2000]{Primary 35B40, 35K15; Secondary 35K05,35K55}
%
%\date{XXXXXXXXXXXXXX}
%%\dedicatory{This paper is dedicated to our advisors.}
%\keywords{large time behavior, heat equation, potential, semilinear heat equation}
%%%%%%%%%%%%%%%%%%%%%%%%%%%%%%%%%%%%%%
%%%%%%%%%%%%%%%%%%%%%%%%%%%%%%%%%%%%%%
\begin{abstract}
In this paper we consider an overdetermined Cauchy problem for the heat equation.
We prove that if the problem has a non-trivial non-negative solution with a certain sequence of similar level sets,
then the solution must be radially symmetric.
\end{abstract}

%%
%{\bf Key words.}
%XXX

%%
%{\bf AMS subject classifications.}
%XXX
%%%%%%%%%%%%%%%%%%%%%%%%%%%%%%%%%%%%%
%%%%%%%%%%%%%%%%%%%%%%%%%%%%%%%%%%%%%
\section{Introduction}
%%%%%%%%%%%%%%%%%%%%%%%%%%%%%%%%%%%%%
%%%%%%%%%%%%%%%%%%%%%%%%%%%%%%%%%%%%%
Consider the unique bounded solution $u=u(x,t)$ of the Cauchy problem for the heat equation:
\begin{equation}
\label{eq:1.1}
 \partial_t u=\Delta u\quad
 \mbox{in}\quad {\mathbb R}^N\times(0,+\infty),\qquad
u=g\ge 0
\quad\mbox{on}\quad {\mathbb R}^N\times\{0\},
\end{equation}
where $N\ge 1$ and
$g$ is a non-trivial bounded non-negative function.
For the initial data $g$,
we denote by $G_0$ the support of $g$, namely $G_0={\rm spt(}g{\rm)}$.
It is well known that
if $g$ is radially symmetric,
then the solution $u$ of \eqref{eq:1.1} must be radially symmetric.

The overdetermined problems, which determine the shape of solutions by using some additional  information of solutions, 
are interesting ones in the study of qualitative properties of solutions of partial differential equations.
In \cite[Corollary~3.2, p.4829]{MPS},
problem \eqref{eq:1.1}, 
where $g$ is replaced by a characteristic function of a bounded open set,
is considered, and
it is shown that if there exists a non-empty stationary isothermic surface of $u$,
then $u$ must be radially symmetric.
In this paper we consider another type of overdetermination.
Precisely
we consider the Cauchy problem \eqref{eq:1.1} which has a solution with
a certain sequence of similar level sets,
and prove the following.
\begin{theorem}
\label{Theorem:1.1}
Let $N\ge 1$ and let $G_0$ be a compact set.
Suppose that there exists a bounded domain $\Omega$ in $\mathbb R^{N}$ with $C^1$ boundary $\partial \Omega$ and satisfying $G_0 \cup \{0\} \subset \Omega$ such that the solution $u$ of \eqref{eq:1.1}
satisfies the following condition:
$$
{\rm(}C{\rm)}
\quad
\left\{
\begin{array}{l}
\displaystyle{
\mbox{
there exist two sequences of positive numbers $\{t_n\}_{n=1}^\infty$ and $\{a_n\}_{n=1}^\infty$ such that}}
\vspace{7pt}\\
\displaystyle{\,\,\,
\mbox{$t_n\uparrow\infty$ as $n\uparrow\infty$ and
$u(t_nx,t_n)=a_n$ for all $n\in{\mathbb N}$ and $x\in\partial \Omega$.}}
\end{array}
\right.
$$
Then $u$ must be radially symmetric with respect to the origin. 
\end{theorem}
\begin{remark}
\label{re:2}
The overdetermination by only one stationary level surface of  solutions of parabolic problems for diffusion equations has been considered since the paper {\rm\cite{MS1}} appeared.
For instance, in {\rm\cite{MS2}} Magnanini and the second author of the present paper considered 
the initial-boundary value problem, where the initial value equals zero and the boundary value equals $1$, for some nonlinear diffusion equation over a domain $\Omega$ with bounded boundary  $\partial\Omega$, and also the Cauchy problem where the initial data is a characteristic function of  the set $\mathbb R^{N}\setminus\Omega$. Then, 
they proved that if
a solution $u$ has  
a surface $\Gamma\subset\Omega$ of codimension $1$ such that,
for some function $a$ : $(0,T)\to{\mathbb R}$, $u(x,t)=a(t)$ for every $(x,t)=\Gamma\times(0,T)$,
where $\Gamma$ is a so-called stationary level surface of $u$,
then $\partial\Omega$ must be a sphere.
See also {\rm\cite{CMS, MPS, MS1, MS3}}.
\end{remark}

The proof of Theorem~\ref{Theorem:1.1} consists of two steps.
In the first step, 
by using condition $(C)$ and the monotonicity of solutions on some exterior domain
(see \eqref{eq:2.3})
we see that $\partial \Omega$ is a sphere with center the origin (see Lemma~\ref{Lemma:2.1} and Proposition \ref{Proposition:2.1}),
and we prove the radial symmetry of the solution in the second step.

Our argument in the first step is also applicable to some overdetermined elliptic boundary value problems over exterior domains in $\mathbb R^N$.

As an example,
we consider the following boundary value problem for some fully nonlinear elliptic equation.
Let $u\in C^1(\mathcal D) \cap C^0(\overline{\mathcal D})$ be the unique viscosity solution of 
\begin{equation}
\label{eq:1.5}
\left\{
\begin{array}{ll}
\displaystyle{F(D^2u,Du,u)=0} & \mbox{in}\quad {\mathcal D},
\vspace{3pt}\\
u>0
& \mbox{in}\quad {\mathcal D},
\vspace{3pt}\\
u=1
& \mbox{on}\quad {\partial \Omega},
\vspace{3pt}\\
u\to 0
& \mbox{as}\quad |x|\to\infty,
\end{array}
\right.
\end{equation}
where $N\ge 2$,
$\Omega\subset{\mathbb R}^N$ is a bounded domain satisfying $0\in\Omega$ with $C^1$ boundary $\partial\Omega$,
and ${\mathcal D}={\mathbb R}^N\setminus {\overline \Omega}$ is also a domain.
Here the nonlinearity $F$ satisfies the following:
\begin{itemize}
\item[(H1)] 
(Regularity) 
$F$ is a continuous function defined on ${\mathcal S}^N({\mathbb R})\times{\mathbb R}^N\times{\mathbb R}$, 
where ${\mathcal S}^N$ denotes the space of $N\times N$ symmetric (real) matrices.
Furthermore,
for any $R>0$
there exists a positive constant $C_1$ such that
$$
|F(M,p,u_1)-F(L,q,u_2)|\le C_1\{|M-L|+|p-q|+|u_1-u_2|\}
$$
for all $M,L\in{\mathcal S}^N(\mathbb R)$, $p,q\in{\mathbb R}^N$, and $u_1,u_2\in[-R,R]$.
\vspace{3pt}
\item[(H2)]
(Ellipticity)
There exists a constant $C_2>0$ such that
$$
F(M+L,p,u)-F(M,p,u)\ge C_2\mbox{Tr}(L)
$$
for all $M,L\in{\mathcal S}^N(\mathbb R)$ with $L\ge 0$, $p\in{\mathbb R}^N$, and $u\in{\mathbb R}$.
\vspace{3pt}
\item[(H3)]
(Symmetry)
For any $M\in{\mathcal S}^N(\mathbb R)$,  $A \in {\mathcal O}^N(\mathbb R)$, $p\in{\mathbb R}^N$, and $u\in(0,\infty)$,
$$
F(M,p,u)=F({}^tAMA, {}^tAp,u),
$$
where ${\mathcal O}^N(\mathbb R)$ denotes the set of $N$-dimensional orthogonal matrices and ${}^tA$ denotes the transpose of $A \in {\mathcal O}^N(\mathbb R)$.
\vspace{3pt}
\item[(H4)](Homogeneity)
There exists some constant $\beta < 0$ such that,
for any $\mu > 1$, the function
\begin{equation}
\label{eq:1.6}
u_\mu(x):=\mu^\beta u(\mu^{-1} x)
\end{equation}
is a solution of problem \eqref{eq:1.5}, where $\Omega$ is replaced  by
$$
\Omega_\mu=\{\mu x\in{\mathbb R}^N : x\in \Omega\}
$$
and where the boundary condition $u=1$ on $\partial\Omega$ is replaced by
$$
u_\mu(x) =  \mu^\beta\ \mbox{ on }\ \partial\Omega_\mu.
$$
\end{itemize}
Then the following holds.
\begin{theorem}
\label{Theorem:1.2}
Suppose that $F$ satisfies {\rm (H1)--(H4)} and moreover $F$ is nonincreasing in $u>0$.
Let $u\in C^1(\mathcal D) \cap C^0(\overline{\mathcal D})$ be a viscosity solution of \eqref{eq:1.5}.
Assume that
there exists a constant $\lambda>1$ such that $\overline{\Omega} \subset \Omega_\lambda$ and
\begin{equation}
\label{eq:1.6-1}
u(x)=\lambda^\beta\quad\mbox{for all}\quad x\in\partial \Omega_\lambda.
\end{equation}
Then $\partial\Omega$ is a sphere with center the origin and $u$ must be radially symmetric.
\end{theorem}
\begin{remark}
\label{re:1}
Such overdetermination by only one level surface of  solutions of elliptic boundary value problems over bounded domains has been considered. 

In {\rm\cite{ES}} Enache and the second author of the present paper studied the overdetermination by only one level surface which is similar to the boundary. To be precise, they considered
some boundary value problem for some fully nonlinear elliptic problem in a bounded domain $\Omega$, 
and proved that, if there exist  constants $\lambda\in(0,1)$ and $\alpha$ such that
$u(x)=\alpha$ on $\partial {\Omega}_\lambda$, 
then $\Omega$ must be the interior of an $N$-dimensional ellipsoid.
See {\rm\cite[{\it Theorem} 2.1]{ES}}.

In {\rm\cite{CMS, S}}, the overdetermination by only one level surface which is parallel to the boundary was considered. By applying the method of moving planes directly to the problems as in {\rm\cite[Proof of Theorem 1.2, pp. 941--942]{MS2}}, the authors proved that the underlying domain must be a ball.

\end{remark}

The paper is organized as follows. 
In Section~2 we give some preliminary proposition, and prove the key lemma of this paper.
By using this lemma and the proposition we prove the main theorems in Section~3.
In Section~4 we give two remarks on condition $(C)$ concerning a sequence of similar level sets.

%%%%%%%%%%%%%%%%%%%%%%%%%%%%%%%%%%%%%
%%%%%%%%%%%%%%%%%%%%%%%%%%%%%%%%%%%%%
\section{Preliminaries}
%%%%%%%%%%%%%%%%%%%%%%%%%%%%%%%%%%%%%
%%%%%%%%%%%%%%%%%%%%%%%%%%%%%%%%%%%%%
We prepare several notations.
For each $r > 0$ and $z \in \mathbb R^{N}$, denote by $B_{r}(z)$ the open ball in $\mathbb R^{N}$ with radius $r$ and center $z$. 
%For each $l\in\partial B_1(0)$,
%let $\Pi_l \subset \mathbb R^{N}$ be the hyperplane with normal $l$ and $0\in\Pi_l$, that is, 
%$$
%\Pi_{l}= \{ x \in \mathbb R^{N} : x\cdot l = 0 \}.
%$$
 For each $C^{1}$ domain $\Omega\subset{\mathbb R}^N$, at $p\in\partial\Omega$, $T_p(\partial\Omega)$ and $\nu(p)$ denote the tangent space of $\partial\Omega$ and the outer unit normal vector to $\partial\Omega$, respectively.

We first prove the following proposition:
\begin{Proposition}
\label{Proposition:2.1}
Let $N \ge 2$ and let $\Omega\subset{\mathbb R}^N$ be a bounded $C^{1}$ domain containing the origin. 
If every point vector $p \in \partial\Omega$ is parallel to the outer unit normal vector $\nu(p)$ to $\partial\Omega$,
%\begin{equation}
%\label{eq:2.1}
%\{p\in\partial\Omega\,|\, T_p(\partial\Omega)\ni l\} = \Pi_l\cap\partial\Omega\ \ \mbox{ for every } l\in\partial B_1(0),
%\end{equation}
then $\partial\Omega$ must be a sphere with center the origin.
\end{Proposition}
{\bf Proof.}
Since $\Omega$ is a bounded $C^1$ domain,
$\partial \Omega$ has finitely many connected components,
and each component is a $C^1$ closed hypersurface embedded in ${\mathbb R}^N$. 
Let $\Gamma$ be a component of $\partial\Omega$. 
Then, for any $p\in\Gamma$,
we have
\begin{equation}
\label{eq:2.2}
p\perp T_p(\Gamma).
\end{equation}

Let $p=p(t)$ be a regular curve on $\Gamma$.
Then, by \eqref{eq:2.2} we obtain
$$
p(t)\perp \frac{d}{dt}p(t)\quad\mbox{for all $t$},
$$
namely
$$
 \frac{d}{dt}(|p(t)|^2)=0\quad\mbox{for all $t$}.
$$
Therefore we see that there exists a positive constant $C$ such that $|p(t)|=C$ for all $t$.
This implies that $\Gamma$ is a sphere with center the origin.
Therefore,
since $\Omega$ is a domain containing the origin,
we see that $\partial\Omega$ must be a sphere with center the origin, 
and Proposition~\ref{Proposition:2.1} follows.
$\Box$
\vspace{3pt}

Next we prove the key lemma for the proofs of Theorems~\ref{Theorem:1.1} and \ref{Theorem:1.2}.
\begin{lemma}
\label{Lemma:2.1}
Let $N \ge 2$ and let $u \in C^{1}(\mathbb R^{N}\times (0, \infty))$, and let $\Omega$ be a bounded domain in $\mathbb R^{N}$ containing the origin.
Suppose that there exists a half-space $H$ of $\mathbb R^{N}$ including $\overline{\Omega}$ such that
\begin{equation}
\label{eq:2.3}
\frac{\partial u}{\partial l}<0\qquad\mbox{in}\quad\left({\mathbb R}^N\setminus H\right)\times (0,\infty),
\end{equation}
where $l$ is the outer unit normal vector to $\partial H$ and suppose the following condition:
$$
{\rm(}C{\rm)}
\quad
\left\{
\begin{array}{l}
\displaystyle{
\mbox{
there exist two sequences of positive numbers $\{t_n\}_{n=1}^\infty$ and $\{a_n\}_{n=1}^\infty$ such that }}
\vspace{7pt}\\
\displaystyle{\,\,\,
\mbox{$t_n\uparrow\infty$ as $n\uparrow\infty$, and
$u(t_nx, t_{n})=a_n$ for all $n\in{\mathbb N}$ and $x\in\partial \Omega$.}}
\end{array}
\right.
$$
If $p \in \partial\Omega$ and $l \in T_{p}(\partial\Omega)$, then $p\cdot l \le 0$.
\end{lemma}
{\bf Proof.} Without loss of generality, set $l=(1,0, \dots, 0)$. Then, since $0 \in \Omega$ and $\overline{\Omega} \subset H$, there exists a positive constant $\lambda$ satisfying $H =\{x\in{\mathbb R}^N : x_1<\lambda\}$.
Suppose that there exists a point $p \in \partial\Omega$ such that
$$
p\cdot l = p_1>0 \mbox{ and } l \in T_{p}(\partial\Omega).
$$
Hence, by condition $(C)$ we have
\begin{equation}
\label{eq:2.4}
 \frac{\partial u}{\partial x_1}(t_np, t_{n})=0
\end{equation}
for all $n \in \mathbb N$.
Since $p_1>0$ and $t_n\uparrow\infty$ as $n\uparrow\infty$,
by \eqref{eq:2.4} we see that there exists a sufficiently large number $n_*$ such that
$t_{n_*}p_1>\lambda$ and
$$
\frac{\partial u}{\partial x_1}(t_{n_*}p, t_{n})=0,
$$
which contradicts \eqref{eq:2.3}. $\Box$
%Furthermore,
%applying the same argument above we can easy obtain $a_1\ge 0$.
%Thus we have $a_1=0$.
%This together with Proposition~\ref{Proposition:2.1} implies that
%$\partial D$ is a sphere with center the origin,
%and the proof of Lemma \ref{Lemma:2.1} is complete.
%$\Box$
\vspace{3pt}
%%

%%%%%%%%%%%%%%%%%%%%%%%%%%%%%%%%%%%%%
%%%%%%%%%%%%%%%%%%%%%%%%%%%%%%%%%%%%%
\section{Proofs of Theorems \ref{Theorem:1.1} and \ref{Theorem:1.2}}
%%%%%%%%%%%%%%%%%%%%%%%%%%%%%%%%%%%%%
%%%%%%%%%%%%%%%%%%%%%%%%%%%%%%%%%%%%%
The purpose of this section is to prove Theorems \ref{Theorem:1.1} and \ref{Theorem:1.2}.
We first prove Theorem~\ref{Theorem:1.2}.
\vspace{3pt}

\noindent
{\bf Proof of Theorem \ref{Theorem:1.2}.}
By \eqref{eq:1.6}, \eqref{eq:1.6-1} and the uniqueness of the solution of \eqref{eq:1.5}, 
we see that 
$$
u(x) = \lambda^\beta u(\lambda^{-1}x)\ \mbox{ for all } x \in \mathbb R^{N} \setminus \overline {\Omega_\lambda}.
$$
Therefore,  setting
$$
t_n = \lambda^{n}\ \mbox{ and } a_n= \lambda^{\beta(n +1)}\ \mbox{ for every } n \in \mathbb N
$$
yields that the solution $u$ of \eqref{eq:1.5} satisfies  
$$
\mbox{$t_n\uparrow\infty$ as $n\uparrow\infty$, and
$u(t_nx)=a_n$ for all $n\in{\mathbb N}$ and $x\in\partial \Omega_\lambda$}.
$$
On the other hand, by applying an argument similar to that in the proof of  \cite[Theorem~1.3]{DS} with Aleksandrov's reflection principle (see \cite{GNN}, for example), the maximum principle and Hopf's boundary point lemma (see \cite{BD} and \cite[Proposition~2.6]{DS}), 
we see that,  for  each direction $l \in \partial B_1(0)$,  if a half-space $H$ of $\mathbb R^{N}$ includes $\overline{\Omega}$ and has the outer unit normal vector $l$ to $\partial H$, then  the solution of \eqref{eq:1.5} satisfies  that 
$$
\frac {\partial u}{\partial l} < 0\quad\mbox{on}\quad \partial H. 
$$
Moreover, since every half-space including the above half-space satisfies the same conditions,  we notice that \eqref{eq:2.3} of Lemma~\ref{Lemma:2.1} holds true also for the solution of \eqref{eq:1.5}. 
Therefore we can apply Lemma~\ref{Lemma:2.1}  to  every direction $l \in \partial B_1(0)$ and conclude that every point vector $p \in \partial\Omega$ is parallel to the outer unit normal vector $\nu(p)$ to $\partial\Omega$.
Hence by Proposition \ref{Proposition:2.1} $\partial \Omega$ must be a sphere with center the origin.
Therefore, since  for every $A \in {\mathcal O}^N(\mathbb R)$ the function $u(Ax)$ also satisfies \eqref{eq:1.5} by (H3), 
then, by the uniqueness of the solution of \eqref{eq:1.5}, $u(x) \equiv u(Ax)$ and hence $u$ must be radially symmetric.
The proof of Theorem~\ref{Theorem:1.2} is complete.
$\Box$
\vspace{5pt}

Next we prove Theorem~\ref{Theorem:1.1}.
Let $H$ be an arbitrary half-space of $\mathbb R^{N}$ including $G_0$.
Then, by Aleksandrov's reflection principle, the maximum principle
and Hopf's boundary point lemma,
we have
\begin{equation}
\label{eq:3.1}
\frac{\partial u}{\partial l}<0\qquad\mbox{on}\quad\partial H \times(0,\infty),
\end{equation}
where $l$ is the outer unit normal vector to $\partial H$. Moreover $u$ satisfies \eqref{eq:2.3} of Lemma~\ref{Lemma:2.1}.
Therefore, by  condition $(C)$ we can use Lemma~\ref{Lemma:2.1} and hence by Proposition~\ref{Proposition:2.1} $\partial\Omega$ must be a sphere with center the origin for $N \ge 2$.

\par
We first prove Theorem~\ref{Theorem:1.1} with $N=1$.
\vspace{5pt}

%%%%%%%%%%%%%%%%%%%%%%%%%%%%%%%%%%%%%
%%%%%%%%%%%%%%%%%%%%%%%%%%%%%%%%%%%%%
%%%%%%    Proof of Theorem 1.1 for N=1   %%%%%%%%%%%%%%
%%%%%%%%%%%%%%%%%%%%%%%%%%%%%%%%%%%%%
%%%%%%%%%%%%%%%%%%%%%%%%%%%%%%%%%%%%%
\noindent
{\bf Proof of Theorem \ref{Theorem:1.1} for $N=1$.}
Since $0 \in G_0 \subset \Omega$, we can set $\Omega = (a, b)$ for some $a < 0 < b$.
Then it follows from condition $(C)$ that
\begin{equation}
\label{similar level sets}
u(t_na,t_n) = u(t_nb,t_n) (= a_n) \mbox{ for every } n \in \mathbb N.
\end{equation}
Let us show that $a + b=0$. Suppose that $a+b > 0$. Then, since $t_n \uparrow \infty$ as $n \uparrow \infty$,  there exists $m \in \mathbb N$ such that
\begin{equation}
\label{Large m}
\frac {t_m(a+b)}2 > b.
\end{equation}
Consider the function $w = w(x,t)$ defined by
$$
w(x,t) = u(x,t) - u(t_ma+t_mb-x, t).
$$
Then we have from \eqref{Large m} the following:
\begin{eqnarray*}
& \partial_t w = \partial_x^2 w\ &\mbox{ in } \left({t_m(a+b)}/2, +\infty\right)\times (0,+\infty),
\\
&w = 0\ &\mbox{ on } \left\{ {t_m(a+b)}/2\right\} \times  (0,+\infty),
\\
&w \le 0 \mbox{ and } w \not\equiv 0\ &\mbox{ on } \left( {t_m(a+b)}/2, +\infty\right)\times \{ 0\}.
\end{eqnarray*}
Thus it follows from the strong maximum principle that 
$$
w < 0 \ \mbox{ in } \left( {t_m(a+b)}/2, +\infty\right)\times (0,+\infty),
$$
which contradicts the fact that $w(t_mb,t_m) = 0$ because of \eqref{similar level sets}.
Therefore, we conclude that  $a+b \le 0$. By the same argument we also conclude that $a+b \ge 0$.

Here we can put $\Omega = (-b, b)$.
For $(x,t)\in[0,\infty)\times[0,\infty)$,
consider the functions $v=v(x,t)$ and $v_0 = v_0(x)$ defined by
$$
v(x,t)=u(x,t)-u(-x,t)\ \mbox{ and } v_0(x) = g(x) -g(-x).
$$
It suffices to prove 
\begin{equation}
\label{eq:3.2}
v_0(x)=0\qquad \mbox{for almost every}\quad x\in[0,\infty).
\end{equation}
Indeed, if \eqref{eq:3.2} holds,
then $g(x)=g(|x|)$ for almost every $x\in{\mathbb R}$.
This together with the uniqueness of the solution yields the conclusion of Theorem \ref{Theorem:1.1} with $N=1$.
 
Since $G_0 \subset \Omega = (-b,b)$,
we see that  ${\rm spt}(v_0)\subset(-b, b)$.
This implies that $v$ satisfies
\begin{equation}
\label{eq:3.3}
v(x,t)=(4\pi t)^{-\frac{1}{2}}\int_{-b}^be^{-\frac{(x-y)^2}{4t}}v_0(y)\,dy
\end{equation}
for all $(x,t)\in{\mathbb R}\times(0, +\infty)$.
Furthermore, by \eqref{similar level sets} with $a = -b$, we see that
$$
v(t_nb,t_n)=0  \mbox{ for every } n \in \mathbb N.
$$
This together with \eqref{eq:3.3} yields
\begin{equation}
\label{eq:3.4}
\int_{-b}^bv_0(y)e^{\frac{by}{2}}e^{-\frac{y^2}{4t_n}}\,dy=0 \mbox{ for every } n \in \mathbb N.
\end{equation}
Put
\begin{equation}
\label{eq:3.5}
w_0(y)=v_0(y)e^{\frac{by}{2}}.
\end{equation}
Then, by \eqref{eq:3.4} we have
$$
\int_0^{b^2}\frac{w_0(\sqrt{s})+w_0(-\sqrt{s})}{2\sqrt{s}}e^{-\frac{s}{4t_n}}\,ds=0 \mbox{ for every } n \in \mathbb N.
$$
Since $t_n\to\infty$ as $n\to\infty$,
by the analyticity of the exponential function
we obtain
$$
\int_0^{b^2}\frac{w_0(\sqrt{s})+w_0(-\sqrt{s})}{2\sqrt{s}}e^{-\lambda s}\,ds=0,\qquad \lambda\in\mathbb R.
$$
This together with the injectivity of the Laplace transform yields
\begin{equation}
\label{eq:3.6}
w_0(\sqrt{s})+w_0(-\sqrt{s})=0\ \mbox{ for almost every } s > 0.
\end{equation}
By \eqref{eq:3.5} and \eqref{eq:3.6} we have
$$
v_0(\sqrt{s})e^{\frac{b\sqrt{s}}{2}}+v_0(-\sqrt{s})e^{-\frac{b\sqrt{s}}{2}}=0 \ \mbox{ for almost every } s > 0.
$$
This implies that $v_0(\sqrt{s})=0$ for almost every $s>0$.
Thus we have \eqref{eq:3.2}, and Theorem~\ref{Theorem:1.1} with $N=1$ follows.
$\Box$
\vspace{3pt}

Next we prove Theorem~\ref{Theorem:1.1} for the case $N\ge 2$.
Before beginning the proof
we recall the following lemma,
which follows from the Funk-Hecke formula (see \cite[Theorem~2.22, p. 36]{AH} or \cite[Theorem 6, p. 20]{M}) and Rodrigues' formula
(see \cite[Theorem~2.23, p. 37]{AH} or \cite[Theorem 5, p. 17]{M}).
%%%%%%%%%%%%%%%%%%%%%%%%%%%%%%%%%%%%%
%%%%%%%%%%%%%%%%%%%%%%%%%%%%%%%%%%%%%
%%%%%%    Lemma 3.1 begins             %%%%%%%%%%%%%%
%%%%%%%%%%%%%%%%%%%%%%%%%%%%%%%%%%%%%
%%%%%%%%%%%%%%%%%%%%%%%%%%%%%%%%%%%%%
\begin{lemma}
\label{Lemma:3.1}
Let $L \not= 0$ be a real constant. For $f \in L^{2}(S^{N-1})$, set
$$
{\mathcal L f}(\omega)=\int_{S^{N-1}}e^{L\alpha\cdot\omega}f(\alpha)\,d\sigma(\alpha) \ \mbox{ for every } \omega \in S^{N-1},
$$
where $d\sigma(\alpha)$ denotes the area element of the $(N-1)$-dimensional unit sphere $S^{N-1}$ in $\mathbb R^{N}$.
Then the set $\{ {\mathcal L}f : f \in  L^{2}(S^{N-1}) \}$  is dense in $L^2(S^{N-1})$.
\end{lemma}
{\bf Proof.} Let $p=p(x)$ be an arbitrary harmonic homogeneous polynomial of degree $k \ge 0$ in $\mathbb R^N$. 
Then it follows from  the Funk-Hecke formula that
$$
{\mathcal L p}(\omega)= \lambda p(\omega) \ \mbox{ for every } \omega \in S^{N-1} \ \mbox{ and } \lambda 
= |S^{N-2}|\int_{-1}^1e^{Lt} P_k(t) (1-t^2)^{\frac{N-3}2}\,dt,
$$
where $|S^{N-2}|$ denotes the volume of the $(N-2)$-dimensional unit sphere in $\mathbb R^{N-1}$ 
and $P_k(t)$ denotes the Legendre polynomial of degree $k$ in $\mathbb R^N$. 
Moreover, Rodrigues' formula gives
$$
P_k(t) = (-1)^k\frac{\Gamma\left(\frac{N-1}2\right)}{2^k\Gamma\left(k+\frac{N-1}2\right)}
(1-t^2)^{\frac{3-N}2}\left(\frac d{dt}\right)^k(1-t^2)^{k+\frac{N-3}2}.
$$
Therefore, integrating by parts $k$ times on the definition of the number $\lambda$ yields that
$$
\lambda =  |S^{N-2}| \frac{\Gamma\left(\frac{N-1}2\right)}{2^k\Gamma\left(k+\frac{N-1}2\right)} L^k \int_{-1}^1e^{Lt}(1-t^2)^{\frac{N-3}2}\,dt  \not= 0.
$$
This implies that the linear space $\{ {\mathcal L}f : f \in  L^{2}(S^{N-1}) \}$ contains all the spherical harmonics 
because any spherical harmonic is given by restricting a harmonic homogeneous polynomial onto $S^{N-1}$. 
Therefore, the conclusion holds true.
$\Box$
\vspace{5pt}

Now we are ready to prove Theorem~\ref{Theorem:1.1} for the case $N\ge 2$.
\vspace{5pt}
%%%%%%%%%%%%%%%%%%%%%%%%%%%%%%%%%%%%%
%%%%%%%%%%%%%%%%%%%%%%%%%%%%%%%%%%%%%
%%%%%%    Proof of Theorem 1.1 for N \ge 2   %%%%%%%%%%%%
%%%%%%%%%%%%%%%%%%%%%%%%%%%%%%%%%%%%%
%%%%%%%%%%%%%%%%%%%%%%%%%%%%%%%%%%%%%

\noindent
{\bf Proof of Theorem \ref{Theorem:1.1} for $N\ge 2$.}
First of all, since we already know that $\partial\Omega$ is a sphere with center the origin, there exists a constant $R>0$ such that
$\Omega = B_R(0)$.
Take $A\in\mathcal O^N(\mathbb R)$ arbitrarily.
Then, it follows from $(C)$ that
\begin{equation}
\label{eq:3.7}
u(t_nx,t_n)-u(t_nAx,t_n)=0,\qquad x\in\partial B_R(0),
\end{equation}
for all $n \in \mathbb N$.
Since $G_0\subset B_R(0)$, by \eqref{eq:3.7} we have
$$
\int_{|y|\le R}e^{\frac{x\cdot y}{2}}e^{-\frac{|y|^2}{4t_n}}(g(y)-g(Ay))\,dy=0, \qquad x \in \partial B_R(0), 
$$
for all $n \in \mathbb N$.
This, together with the analyticity of the exponential function and the fact that
$t_n\to\infty$ as $n\to\infty$, implies that
\begin{equation}
\label{eq:3.8}
\int_{|y|\le R}e^{\frac{x\cdot y}{2}}e^{-s|y|^2}(g(y)-g(Ay))\,dy=0, \qquad x\in\partial B_R(0),
\end{equation}
for all $s\in{}\mathbb R$. 
By setting $x = R \alpha$ for $\alpha\in{S}^{N-1}(= \partial B_1(0))$, we have from
 \eqref{eq:3.8} that
$$
\int_0^Rr^{N-1}e^{-sr^2}\int_{S^{N-1}}e^{\frac{Rr}{2}\alpha\cdot\omega}(g(r\omega)-g(rA\omega))\,d\sigma(\omega)\,dr=0, 
\qquad \alpha\in{S}^{N-1},
$$
for all $s\in{\mathbb R}$.
This together with the injectivity of the Laplace transform yields
\begin{equation}
\label{eq:3.9}
\int_{S^{N-1}}e^{\frac{Rr}{2}\alpha\cdot\omega}(g(r\omega)-g(rA\omega))\,d\sigma(\omega)=0, \qquad \alpha\in{S}^{N-1},
\end{equation}
for almost every $r> 0$. 
Let $f \in L^{2}(S^{N-1})$.
Then, by \eqref{eq:3.9} we obtain
$$
\int_{S^{N-1}}\int_{S^{N-1}}e^{\frac{Rr}{2}\alpha\cdot\omega}(g(r\omega)-g(rA\omega))f(\alpha)\,d\sigma(\omega)\,d\sigma(\alpha)=0
$$
for almost every $r> 0$.
Thus, setting $L = \frac{Rr}{2}$ for almost every fixed $r> 0$
 in Lemma~\ref{Lemma:3.1} yields that for almost every fixed $r> 0$
$$
\int_{S^{N-1}}{\mathcal L}f(\omega)(g(r\omega)-g(rA\omega))\,d\sigma(\omega)=0\ \mbox{ for every } f \in  L^{2}(S^{N-1}).
$$
Then, by Lemma~\ref{Lemma:3.1} we have that, for almost every fixed $r> 0$,
$$
g(r\omega)=g(rA\omega)\ \mbox{ for almost every } \omega \in S^{N-1}.
$$
This yields that $g(x) = g(Ax)$ for almost every $x \in \mathbb R^{N}$ and hence by the uniqueness of the solution we have that 
$$
u(x,t) = u(Ax,t)\ \mbox{ for all } (x,t) \in \mathbb R^{N} \times (0,\infty).
$$
Therefore, since $A \in \mathcal O^N(\mathbb R)$ is arbitrary, we see that $u$ is  radially symmetric with respect to the origin. This completes the proof of Theorem \ref{Theorem:1.1} for $N\ge 2$.
$\Box$
%%

%%%%%%%%%%%%%%%%%%%%%%%%%%%%%%%%%%%%%
%%%%%%%%%%%%%%%%%%%%%%%%%%%%%%%%%%%%%
\section{Remarks on condition $(C)$}
%%%%%%%%%%%%%%%%%%%%%%%%%%%%%%%%%%%%%
%%%%%%%%%%%%%%%%%%%%%%%%%%%%%%%%%%%%%

In this last section we give two remarks on condition $(C)$ concerning a sequence of similar level sets.
First one means that 
the order of dependence for the sequence $\{t_n\}_{n=1}^\infty$ with respect to spatial variables in condition $(C)$ is important.
Second one says that 
if the solution $u$ has similar level sets continuously with time,
then we can easily carry out the second step of the proof of Theorem~\ref{Theorem:1.1}.
\begin{remark}
\label{Remark:4.1}
There exists a solution of \eqref{eq:1.1} with $N=3$ which is not radially symmetric
even if it has a sequence of similar level sets. 

\vskip 2ex
Let $a$ be a positive constant and
$v_0 \in C_0^\infty(\mathbb R)$ be a nonnegative even function satisfying
\begin{equation}
\label{eq:4.0}
{\rm spt(}v_0{\rm)}=[-a,a]\quad\mbox{ and }\quad
v_0'<0\quad\mbox{if}\quad s\in(0,a).
\end{equation}
Put
\begin{equation}
\label{eq:4.1}
v(s,t)=(4\pi t)^{-\frac{1}{2}}\int_{-a}^ae^{-\frac{(s-\mu)^2}{4t}}v_0(\mu)\,d\mu,
\qquad
w(s,t)=\frac{\partial^3}{\partial s^3}v(s,t),
\end{equation}
for $(s,t)\in{\mathbb R}\times(0,\infty)$.
Then, since $v$ is a even function in $s$,
$w$ is a odd function in $s$.
Furthermore, we set $r=|x|$ for $x\in{\mathbb R}^3$ and put
\begin{equation}
\label{eq:4.2}
f(r,t)=\frac{\partial}{\partial r}\left(\frac{w(r,t)}{r}\right)=\left(r\frac{\partial w}{\partial r}-w\right)r^{-2}.
\end{equation}
Then we have the following lemma:
\begin{lemma}
\label{Lemma:4.1}
Let $f$ be the function given in \eqref{eq:4.2}.
Then there exists a positive function $r(t)$ for $t > 0$ such that
$$
f(r(t),t)=0,\qquad r(t)=O(t^{\frac{1}{2}})\quad\mbox{as}\quad t\to\infty.
$$
\end{lemma}
{\bf Proof.}
By \eqref{eq:4.0} and \eqref{eq:4.1}
we can take three positive functions $r_2(t), \ r_3(t)$ and $r_4(t)$ for $t > 0$  such that $r_2(t) < r_3(t)<r_4(t)$ and
\begin{eqnarray}
\label{eq:a}
&&
\frac{\partial^3 }{\partial s^3}v(r_3(t),t)=0,\qquad \frac{\partial^3 }{\partial s^3}v(s,t)<0\quad\mbox{ if }\quad s>r_3(t),
\\
\label{eq:b}
&&
\frac{\partial^4 }{\partial s^4}v(r_2(t),t)=\frac{\partial^4 }{\partial s^4}v(r_4(t),t)=0,
\\
\label{eq:c}
&&
 \frac{\partial^4 }{\partial s^4}v(s,t)<0\quad\mbox{ if }\quad r_2(t) < r < r_4(t),\ \mbox{ and }
 \frac{\partial^4 }{\partial s^4}v(s,t)>0\quad\mbox{for}\quad s>r_4(t),
\end{eqnarray}
for all $t>0$.
Put
$$
 h(s,t)=s^2 f(s,t) =s\frac{\partial^4 v}{\partial s^4}-\frac{\partial^3 v}{\partial s^3}
$$
Then, since $h(s,t)>0$ for $s\ge r_4(t)$ and $h(s,t)<0$ for $s\in[r_2(t),r_3(t)]$,
by applying the intermediate value theorem,
we can take a positive function $r(t)$ for $t > 0$  such that
\begin{equation}
\label{eq:4.4}
f(r(t),t)=0,\qquad r_3(t)<r(t)<r_4(t).
\end{equation}
On the other hands,
by \eqref{eq:4.1} we obtain
\begin{eqnarray}
\label{eq:4.5}
&&
(4\pi t)^{\frac{1}{2}}\frac{\partial^3 v}{\partial s^3}
=-\frac{1}{8t^3}\int_{-a}^a(s-\mu)\{-6t+(s-\mu)^2\}e^{-\frac{(s-\mu)^2}{4t}}v_0(\mu)\,d\mu,
\\
&&
(4\pi t)^{\frac{1}{2}}\frac{\partial^4 v}{\partial s^4}
=\frac{1}{16t^4}\int_{-a}^a\{12t^2-12t(s-\mu)^2+(s-\mu)^4\}e^{-\frac{(s-\mu)^2}{4t}}v_0(\mu)\,d\mu
\nonumber\\
\label{eq:4.6}
&&\qquad\qquad\quad
>\frac{1}{16t^4}\int_{-a}^a(s-\mu)^2\{(s-\mu)^2-12t\}e^{-\frac{(s-\mu)^2}{4t}}v_0(\mu)\,d\mu.
\end{eqnarray}
for all $(s,t)\in{\mathbb R}\times(0,\infty)$.
By \eqref{eq:4.5} we see that
\begin{eqnarray*}
&&
\frac{\partial^3v}{\partial s^3}<0\quad\mbox{for}\quad s>a+\sqrt{6t},
\\
&&
\frac{\partial^3v}{\partial s^3}>0\quad\mbox{for}\quad a<s<\sqrt{5t}-a\quad\mbox{with}\quad t>\frac{4a^2}{5}.
\end{eqnarray*}
These together with \eqref{eq:a} yield
\begin{equation}
\label{eq:4.7}
\sqrt{5t}-a\le r_3(t)\le a+\sqrt{6t}\quad\mbox{for}\quad t>\frac{4a^2}{5}.
\end{equation}
Furthermore, by \eqref{eq:4.6} we have
$$
\frac{\partial^4v}{\partial s^4}>0\quad\mbox{for}\quad s>a+\sqrt{12t}.
$$
This together with \eqref{eq:b} implies that
\begin{equation}
\label{eq:4.8}
r_4(t)\le a+\sqrt{12t}.
\end{equation}
Combining \eqref{eq:4.4}, \eqref{eq:4.7} and \eqref{eq:4.8} yields that
$$
\sqrt{5t}-a\le r_3(t)<r(t)<r_4(t)\le a+\sqrt{12t}
$$
for all $t>(4a^2)/5$.
This implies that $r(t)=O(t^{\frac{1}{2}})$ as $t\to\infty$;
thus Lemma~{\rm\ref{Lemma:4.1}} follows.
$\Box$
\vspace{8pt}

On the other hand, by \eqref{eq:4.0} and \eqref{eq:4.2}
we can take a nonnegative radially symmetric function $\psi(|x|)\in C_0^\infty({\mathbb R}^3)$ such that
$$
\psi(|x|)+f(|x|,0)\frac{x_1}{|x|}\ge 0
$$
for all $x\in{\mathbb R}^3$.
We set
\begin{equation}
\label{eq:4.3}
u(x,t)=(4\pi t)^{-\frac{3}{2}}\int_{{\mathbb R}^3}e^{-\frac{|x-y|^2}{4t}}\psi(|y|)\,dy+f(|x|,t)\frac{x_1}{|x|}
:=u_{rad}(|x|,t)+f(|x|,t)\frac{x_1}{|x|}.
\end{equation}
Then the function $u$ is a solution of \eqref{eq:1.1} with $N =3$ and
$$
g(x)=\psi(|x|)+f(|x|,0)\frac{x_1}{|x|}.
$$
By Lemma~{\rm\ref{Lemma:4.1}} and \eqref{eq:4.3} we see that, if $|x|=r(t)$, then
there exists a function $c(t)$ such that
$$
u(x,t)=u_{rad}(r(t),t)=c(t).
$$
This implies that
the solution $u$ is not radially symmetric even if it has a sequence of similar level sets.
\end{remark}

\begin{remark}
\label{Remark:4.2}
Instead of $(C)$,
suppose that there exists a function $a=a(t)$ for $t>0$ such that
$$
u((1+t)x,t)=a(t)
$$
for all $(x,t)\in\partial \Omega\times(0,\infty)$.
Then we can use the maximum principle and the unique continuation theorem 
and get the same conclusion of Theorem $\ref{Theorem:1.1}$.

\vskip 2ex
Indeed, as in the proof of Theorem~$\ref{Theorem:1.1}$,  
by Aleksandrov's reflection principle, the maximum principle, Hopf's boundary point lemma, 
Proposition~$\ref{Proposition:2.1}$ and Lemma~$\ref{Lemma:2.1}$, 
we see that $\partial\Omega$ must be a sphere with center the origin for $N \ge 2$.  
Say $\partial\Omega = \partial B_R(0)$ for some $R > 0$. 
Take $A \in \mathcal O^N(\mathbb R)$ arbitrarily. 
Then
$$
u((1+t)x,t) - u((1+t)Ax,t) = 0,\qquad x \in \partial B_R(0),\ t > 0.
$$
Since $u((1+0)x,0) - u((1+0)Ax,0) = 0$ if $x \not\in B_R(0)$, by the maximum principle we get
$$
u(x,t) - u(Ax,t) = 0 \quad \mbox{ if } \ |x| > (1+t)R.
$$
Hence it follows from the unique continuation theorem $($see {\rm\cite{EF}}$)$ 
that $u(x,t) - u(Ax,t)$ equals zero identically, which gives the conclusion.
\end{remark}
\noindent
{\bf Acknowledgements.}
The first and second authors were supported by the Grants-in-Aid for Young Scientists (B) (No.~24740107)
and for Challenging Exploratory Research (No.~25610024) from Japan Society for the Promotion of Science, respectively. The authors would like to thank the anonymous referee for his/her some
valuable suggestions to improve the presentation and clarity in several points.

%%%%%%%%%%%%%%%%%%%%%%%%%%%%%%%%%%%%%%
%%%%%%%%%%%%    references    %%%%%%%%%%%%%%%%%%
%%%%%%%%%%%%%%%%%%%%%%%%%%%%%%%%%%%%%%
\bibliographystyle{amsplain}

\begin{thebibliography}{10}


\bibitem{AH}
K. Atkinson and W. Han, 
Spherical Harmonics and Approximations on the Unit Sphere: An Introduction, Lecture Notes in Math.  {\bf 2044}, Springer-Verlag, 2012.

\bibitem{BD}
M. Bardi and F. Da Lio, 
On the strong maximum principle for fully nonlinear degenerate elliptic equations, 
Arch. Math. {\bf 73} (1999), 276--285.

%\bibitem{BKP}
%M. Bertsch, R. Kersner, and L. A.  Peletier, 
%Positivity versus localization in degenerate diffusion equations, 
%Nonlinear Anal. {\bf 9} (1985), 987--1008.

%\bibitem{BN}
%H. Berestycki and L. Nirenberg, 
%On the method of moving planes and the sliding method, 
%Bol. Soc. Brasil. Mat. (N.S.) {\bf 22} (1991), 1--37.

\bibitem{CMS}
G. Ciraolo, R. Magnanini and S. Sakaguchi,
Symmetry of minimizers with a level surface parallel to the boundary, 
to appear in J. Eur. Math. Soc.

\bibitem{DS}
F. Da Lio and B. Sirakov, 
Symmetry results for viscosity solutions of fully nonlinear uniformly elliptic equations, 
J. Eur. Math. Soc. {\bf 9} (2007), 317--330.

\bibitem{ES}
C. Enache and S. Sakaguchi,
Some fully nonlinear elliptic boundary value problems with ellipsoidal free boundaries,
Math. Nachr. {\bf 284} (2011), 1872--1879.


\bibitem{EF} 
L. Escauriaza and F. J. Fern\'{a}ndez,
Unique continuation for parabolic operators, 
Ark. Mat. {\bf 41} (2003), 35--60.

\bibitem{GNN}
B. Gidas, W.-M. Ni and L. Nirenberg, 
Symmetry and related properties via the maximum principle,
Comm. Math. Phys. {\bf 68} (1979), 209--243.

%\bibitem{G}
%H. Groemer, Geometric Applications of Fourier Series and Spherical Harmonics,
%Encyclopedia of Mathematics and its Applications {\bf 61} Cambridge Univ. Press  1996.

%\bibitem{H}
%S. Helgason,
%Groups and geometric analysis. 
%Integral geometry, invariant differential operators, and spherical functions. 
%Mathematical Surveys and Monographs, {\bf 83}. 
%American Mathematical Society, Providence, RI, 2000. 

\bibitem{MPS}
R. Magnanini, J. Prajapat and S. Sakaguchi, 
Stationary isometric surfaces and uniformly dense domains, 
Trans. Amer. Math. Soc. 
{\bf 358} (2006), 4821--4841.

\bibitem{MS1}
R. Magnanini and S. Sakaguchi, 
Matzoh ball soup: Heat conductors with a stationary isothermic surface, 
Ann. of Math.
{\bf 156} (2002), 931--946.

\bibitem{MS2}
R. Magnanini and S. Sakaguchi, 
Nonlinear diffusion with a bounded stationary  level surface, 
Ann. Inst. H. Poincar\'e Anal. Non Lin\'eaire.
{\bf 27} (2010), 937--952.

\bibitem{MS3}
R. Magnanini and S. Sakaguchi, Matzoh ball soup revisited: the boundary regularity issue,  
Math. Methods Appl. Sci. {\bf 36} (2013), 2023--2032.

\bibitem{M}
C. M\"uller, Spherical Harmonics, Lecture Notes in Math.  {\bf 17}, Springer-Verlag, 1966.

\bibitem{S}
H. Shahgholian,
Diversifications of Serrin's and related symmetry problems, 
Complex Var. Elliptic Equ. {\bf 57} (2012), 653--665.


%\bibitem{S}
%R. T. Seeley,
%Spherical harmonics, 
%Amer. Math. Monthly {\bf 73} (1966), 115--121.

\end{thebibliography}

\end{document}